\documentclass[11pt,reqno]{amsart}

\usepackage{amssymb}

\catcode`\@=11

\long\def\@savemarbox#1#2{\global\setbox#1\vtop{\hsize\marginparwidth 
  \@parboxrestore\tiny\raggedright #2}}
\newcommand\lref[1]{\ref{#1}%
\@ifundefined{r@DisplaY #1}{}{ (#1)}}

\newcommand\fakelabel[2]{\@bsphack\if@filesw {\let\thepage\relax
   \newcommand\protect{\noexpand\noexpand\noexpand}%
\xdef\@gtempa{\write\@auxout{\string
      \newlabel{#1}{{#2}{\thepage}}}}}\@gtempa
   \if@nobreak \ifvmode\nobreak\fi\fi\fi\@esphack}

\catcode`\@=12

\def\Empty{}
\newcommand\oplabel[1]{
  \def\OpArg{#1} \ifx \OpArg\Empty {} \else
        \label{#1}
  \fi}
\newtheorem{theoremSt}{Theorem}[section]

\newtheorem{exampleSt}[theoremSt]{Example}
\newtheorem{exerciseSt}[theoremSt]{Exercise}

\newcommand\MakeStEnv[1]{
  \newenvironment{#1}[1]{
  \begin{#1St} \oplabel{##1}%
  \global\def\CrntSt{\thetheoremSt}%
}{ 
  \end{#1St} }
  \newenvironment{#1+}[1]{
  \begin{#1St} \label{##1}%
  \label{DisplaY ##1}%
  \global\def\CrntSt{\thetheoremSt}%
  \def\Labl{##1}\ifx\Labl\Empty{} \else {\em (\Labl)\,}\fi%
}{ 
  \end{#1St} }
}
\MakeStEnv{theorem}
\MakeStEnv{corollary}
\MakeStEnv{proposition}
\MakeStEnv{lemma}
\MakeStEnv{definition}
\MakeStEnv{conjecture}
\MakeStEnv{problem}
\MakeStEnv{question}

\long\def\state#1#2{
\medskip\par\noindent
{\bf #1} 
{\it #2}
\par\medskip
}

\newlength{\saveu}

\newenvironment{pf*}[1]{%
 \begin{proof}[#1]%
}{ 
 \end{proof}
}

\newcommand{\finishproof}[1]{ 
  \def\FPArg{#1}
  \ifx\FPArg\Empty
        \newcommand\FPArg{\CrntSt}  \fi
  \smallbreak\noindent\makebox[\textwidth]{\hfill\fbox{\FPArg}}
  \medbreak\noindent
}

\newcommand\AAA{{\mathcal A}}

\newcommand\CC{{\mathcal C}}

\newcommand\FF{{\mathcal F}}

\newcommand\LL{{\mathcal L}}
\newcommand\MM{{\mathcal M}}

\newcommand\PP{{\mathcal P}}
\newcommand\QQ{{\mathcal Q}}

\newcommand\PMF{{\PP\kern-2pt\MM\FF}}

\newcommand\PML{{\PP\kern-2pt\MM\LL}}

\newcommand\half{{\textstyle{\frac12}}}

\newcommand\ep{\epsilon}

\newcommand\hhat{\widehat}

\newcommand\union{\cup}
\newcommand\intersect{\cap}
\newcommand\bbR{{\mathord{\text{I\kern-2pt R}}}}        
\newcommand\bbH{{\mathord{\text{I\kern-2pt H}}}}        

\newcommand\Z{{\mathbb Z}}
\newcommand\R{{\mathbb R}}
\newcommand\N{{\mathbb N}}

\newcommand\bigrightarrow[1]{\hbox to #1{\rightarrowfill}}
\newcommand\bigleftarrow[1]{\hbox to #1{\leftarrowfill}}

\newcommand\homeo{\cong}
\newcommand\boundary{\partial}
\newcommand\semidir{\mathrel{\hbox{\vrule depth-.03ex height1.1ex\kern-0.15em$\times$}}}

\newcommand{\diam}{\operatorname{diam}}

\numberwithin{equation}{section}

\catcode`\@=11

\def\subsection{\@startsection{subsection}{2}%
  \z@{.5\linespacing\@plus.7\linespacing}{.5em}%
  {\normalfont\bfseries\centering}}

\def\section{\@startsection{section}{1}%
  \z@{.7\linespacing\@plus\linespacing}{.5\linespacing}%
  {\normalfont\large\bfseries\centering}}

\def\subsubsection{\@startsection{subsubsection}{3}%
  \z@{.5\linespacing\@plus.7\linespacing}{-.5em}%
  {\normalfont\bfseries}}

\catcode`\@=12

\newcommand{\collar}{\operatorname{\mathbf{collar}}}

\newcommand{\dist}{\operatorname{dist}}

\newcommand{\fsubd}{\mathrel{{\scriptstyle\searrow}\kern-1ex^d\kern0.5ex}}
\newcommand{\bsubd}{\mathrel{{\scriptstyle\swarrow}\kern-1.6ex^d\kern0.8ex}}
\newcommand{\fsubeq}{\mathrel{\raise-.7ex\hbox{$\overset{\searrow}{=}$}}}
\newcommand{\bsubeq}{\mathrel{\raise-.7ex\hbox{$\overset{\swarrow}{=}$}}}

\newcommand{\base}{\operatorname{base}}

\newcommand{\tsh}[1]{\left\{\kern-.9ex\left\{#1\right\}\kern-.9ex\right\}}
\newcommand{\Tsh}[2]{\tsh{#2}_{#1}}

\theoremstyle{remark} \newtheorem{remark}[theoremSt]{Remark}

\newcommand\AM{{\mathcal{M}^{\omega}}}
\newcommand\AP{{\mathcal{P}^{\omega}}}

\newcommand\AQ{{\mathcal{Q}^{\omega}}}

\newcommand\indhat{{\hhat{\operatorname{ind}}}}
\newcommand\ind{{\operatorname{ind}}}
\newcommand\Ind{{\operatorname{Ind}}}
\newcommand\Dim{{\operatorname{dim}}}
\newcommand\Indhat{{\hhat{\operatorname{Ind}}}}
\newcommand\dimhat{{\hhat{\operatorname{dim}}}}
\def\MCG{\mathcal {MCG}}
\newcommand\tint[2]{#1 \pitchfork #2 \ne \emptyset}
\newcommand\notint[2]{#1 \pitchfork #2 = \emptyset}
\def\ulim{\lim_\omega}

\def\cone{{\rm{Cone}}_{\omega}}
\def\dist{{\rm{dist}}}
\def\co{\colon}

\newcommand\seq[1]{\mbox{\boldmath$#1$}}
\newcommand\subseq[1]{\mbox{\boldmath$\scriptstyle#1$}}

\newcommand\balpha{\seq\alpha}

\begin{document}

\title{Dimension and rank for mapping class groups}
\author{Jason A. Behrstock}
\address{University of Utah}
\thanks{First author supported by NSF grants DMS-0091675 and DMS-0604524.}
\email{jason@math.utah.edu}
\author{Yair N. Minsky}
\address{Yale University}
\thanks{Second author supported by NSF grant DMS-0504019.}
\email{yair.minsky@yale.edu}
\thanks{To appear in \emph{Annals of Mathematics}.}

\dedicatory{Dedicated to the memory of Candida Silveira.}

\begin{abstract}
    We study the large scale geometry of the mapping class group,
    $\MCG$.  Our main result is that for any asymptotic cone of
    $\MCG$, the maximal dimension of locally compact subsets coincides
    with the maximal rank of free abelian subgroups of $\MCG$.  An
    application is a proof of Brock-Farb's Rank
    Conjecture which asserts that $\MCG$ has quasi-flats of
    dimension~$N$ if and only if it has a rank~$N$ free abelian
    subgroup. (Hamenstadt has also given a proof of this conjecture,
    using different methods.) We also compute the maximum dimension of
    quasi-flats in Teichmuller space with the Weil-Petersson metric.
\end{abstract}

\maketitle

\bigskip

The coarse geometric structure of a finitely generated group can be
studied by passage to its {\em asymptotic cone}, which is a space
obtained by a limiting process from sequences of rescalings of the
group. This has played an important role in the quasi-isometric
rigidity results of \cite{DrutuSapir:TreeGraded}, 
\cite{KapovichLeeb:haken}
\cite{KleinerLeeb:buildings}, and others. In this paper we study the
asymptotic cone $\AM(S)$ of the mapping class group of a surface of
finite type. Our main result is 

\state{Dimension Theorem.}{%
The maximal topological dimension of a
locally-compact subset of the asymptotic cone of a mapping class
group is equal to the maximal rank of an
abelian subgroup.}

Note that \cite{BirmanLubotzkyMcCarthy} showed that the maximal rank
of an abelian subgroup of a mapping class group 
of a surface with negative Euler characteristic is $3g-3+p$ where $g$
is the genus and $p$ the number of boundary components. This is also
the number of components of a pants decomposition and hence the
largest rank of a pure Dehn twist subgroup.

As an application we obtain a proof of 
the ``geometric rank conjecture'' for mapping class groups, 
formulated
by Brock and Farb \cite{BrockFarb:curvature}, which states:

\state{Rank Theorem.}{%
The geometric rank of the mapping class group of a surface of finite
type is equal to the maximal rank of an abelian subgroup. 
}

Hamenst\"adt has previously announced a proof of the rank conjecture 
for mapping class groups, which has now appeared in 
\cite{Hamenstadt:TT3}. Her 
proof uses the geometry of train tracks and establishes a
homological version of the dimension theorem. Our methods are quite 
different from hers, and we hope that they will be of independent 
interest.

The geometric rank of a group $G$ is defined as the largest $n$ for which
there exists a quasi-isometric embedding $\Z^n \to G$, also known as
an $n$-dimensional quasi-flat. It was proven in 
\cite{FarbLubotzkyMinsky} that, in the mapping class group, 
maximal rank abelian subgroups are quasi-isometrically 
embedded---thereby giving a lower bound on the geometric rank.
This was known when the Rank Conjecture was formulated, thus the 
conjecture was that the known lower bound for the geometric rank 
is sharp.
The affirmation of this conjecture follows immediately 
from the dimension theorem and the observation that 
a quasi-flat, after passage to the asymptotic cone, becomes a 
bi-Lipschitz-embedded copy of $\R^n$. 

We note that in general the maximum rank of (torsion-free)  
abelian subgroups of a given group does 
not yield either an upper or a lower bound on the geometric rank of 
that group.
For instance, non-solvable Baumslag-Solitar groups have 
geometric rank one \cite{Burillo:rankBS}, but contain  
rank two abelian subgroups. 
To obtain groups with geometric rank one, but no 
subgroup isomorphic to $\Z$, 
one may take any finitely generated infinite torsion group.
The $n$-fold product of such a group with itself has $n$-dimensional
quasi-flats, but no copies of $\Z^n$. 

\medskip

Similar in spirit to the above results, and making use of Brock's 
combinatorial model for the Weil-Petersson metric \cite{Brock:wp},  
we also prove:

\state{Dimension Theorem for Teichm\"{u}ller space.}
{Every locally-compact subset of an asymptotic cone of Teichm\"{u}ller 
space with the Weil-Petersson metric  
has topological dimension at most $\lfloor\frac{3g+p-2}{2}\rfloor$.}

The dimension theorem implies the following, which settles another 
conjecture of Brock--Farb.

\state{Rank Theorem for Teichm\"{u}ller space.}
{The geometric rank of the Weil-Petersson metric on the 
Teichm\"{u}ller space of a surface of finite
type is equal to $\lfloor\frac{3g+p-2}{2}\rfloor$. 
}

This conjecture was made by Brock--Farb after proving this result in 
the case $\lfloor\frac{3g+p-2}{2}\rfloor\leq 1$, by showing that in 
such cases Teichm\"{u}ller space is $\delta$-hyperbolic  
\cite{BrockFarb:curvature}. (Alternate proofs of this result were 
obtained in \cite{Behrstock:asymptotic} and \cite{Aramayona:thesis}.)
We also note that the lower bound on the geometric rank of 
Teichm\"{u}ller space is obtained in \cite{BrockFarb:curvature}.

\subsection*{Outline of the proof}

For basic notation and background see \S\ref{background}.

We will define a family $\PP$ of subsets of $\AM(S)$ with the
following properties:  Each $P\in \PP$ comes equipped with a
bi-Lipschitz homeomorphism to a product $F\times \AAA$, where 
\begin{enumerate}
\item $F$ is an $\R$-tree
\item $\AAA$ is the asymptotic cone of the mapping class group of 
a (possibly disconnected) proper subsurface of $S$.
\end{enumerate}
There will also be a Lipschitz map $\pi_P\co \AM(S) \to F$ such that:
\begin{enumerate}
\item The restriction of $\pi_P$ to $P$ is projection to the first
  factor.
\item $\pi_P$ is locally constant in the complement of $P$.
\end{enumerate}
These properties immediately imply that the subsets $\{t\}\times \AAA$
in $P=F\times\AAA$ separate $\AM(S)$ globally.

The family $\PP$ will also have the property that it {\em separates
  points}, that is: for every $x\ne y$ in $\AM(S)$ there exists
  $P\in\PP$ such that $\pi_P(x) \ne \pi_P(y)$. 

Using  induction, we will be able to show that locally compact subsets
of $\AAA$ have dimension at most $r(S)-1$, where $r(S)$ is the
expected rank for $\AM(S)$. The separation properties
above together with a short lemma in dimension theory then imply that
locally compact subsets of $\AM(S)$ have dimension at most $r(S)$. 

\medskip

Section \ref{background} will detail some background material on
asymptotic cones and on the constructions used in Masur-Minsky
\cite{MasurMinsky:complex1,MasurMinsky:complex2} to study the coarse
structure of the mapping class group. Section \ref{Qregions}
introduces product regions in the group and in its asymptotic cone
which correspond to cosets of curve stabilizers. 

Section \ref{Pregions} introduces
the $\R$-trees $F$, which were initially studied by
Behrstock in \cite{Behrstock:asymptotic}. The regions $P\in\PP$ will
be constructed as subsets of the product regions of Section
\ref{Qregions}, in which one factor is restricted to a subset which is
one of the $\R$-trees. The main technical result of the paper is
Theorem \ref{global projection to F}, which constructs the projection
maps $\pi_P$ and establishes their locally-constant properties. 
An almost immediate consequence is Theorem \ref{separate points},
which gives the family of separating sets whose dimension will be
inductively controlled. 

Section~\ref{dimension} applies Theorem~\ref{separate points} to prove
the Dimension Theorem. 

Section~\ref{pants} applies the same techniques
to prove a similar dimension bound for the asymptotic cone of a space 
known as the pants graph and to deduce a
corresponding geometric rank statement there as well. 
These can be 
translated into results for Teichm\"{u}ller space with its Weil-Petersson
metric, by applying Brock's quasi-isometry \cite{Brock:wp} between the
Weil-Petersson metric and the pants graph.

\subsubsection*{Acknowledgements}
The authors are grateful to Lee Mosher for many insightful 
discussions, and for a simplification to the original proof of 
Theorem~\ref{global projection to F}. We would also like to thank 
Benson Farb for helpful comments on an earlier draft.

\section{Background}
\label{background}

\subsection{Surfaces}
Let $S=S_{g,p}$ be a orientable compact connected surface of genus
$g$ and $p$ boundary components. 
The mapping class group, $\MCG(S)$, is defined to be 
$Homeo^{+}(S)/Homeo_0(S)$, the orientation-preserving 
homeomorphisms up to isotopy. This group is finitely-generated 
\cite{Dehn:Translations, Birman:Braids} and for any finite
generating set one considers the word metric in the usual way 
\cite{Gromov:ICMAddress}, whence yielding a metric space which is 
unique up to quasi-isometry. 

Throughout the remainder, we tacitly exclude the case of the closed torus
$S_{1,0}$.  Nonetheless, the Dimension Theorem does hold in this case
since $\MCG(S_{1,0})$ is virtually free so its asymptotic cones are
all one dimensional and the largest rank of its free abelian subgroups
is one.

Let $r(S)$ denote the largest rank of an abelian subgroup of
$\MCG(S)$ when $S$ has negative Euler characteristic.  
In \cite{BirmanLubotzkyMcCarthy}, it was computed that 
$r(S)=3g-3+p$ 
and it is easily seen that this rank is realized by 
any subgroup generated by Dehn twists on a maximal set of disjoint essential
simple closed curves. Moreover, such subgroups are known to be 
quasi-isometrically embedded by results in \cite{Mosher:automatic}, 
when $S$ has punctures, and by \cite{FarbLubotzkyMinsky} in the general 
case.

For an annulus let $r=1$. For a disconnected subsurface
$W\subset S$, with each 
component homotopically essential and not homotopic into the boundary,
and no two annulus components homotopic to each other, let
$r(W)$ be the sum of $r(W_i)$ over the components of $W$.
We note that $r$ is automatically additive over disjoint unions, and
is monotonic with respect to inclusion.

\subsection{Quasi-isometries}

If $(X_{1},d_{1})$ and $(X_{2},d_{2})$ are metric spaces, 
a map $\phi\co X_{1}\to X_{2}$ is called a 
\emph{$(K,C)$-quasi-isometric embedding} if for each $y,z\in X_{1}$ 
we have: 
\begin{equation}
d_{2}(\phi(y),\phi(z)) \approx_{K,C} d_1(y,z).
\end{equation}
Here the expression $a\approx_{K,C} b$  means
$a/K - C \le b \le Ka+C$. We sometimes suppress $K,C$, writing just
$a\approx b$ when this will not cause confusion.

We call $\phi$ a \emph{quasi-isometry} if, additionally, there 
exists a constant $D\geq 0$ so that each $q\in X_{2}$ satisfies 
$d_{2}(q,\phi(X_{1}))\leq D$, i.e., $\phi$ is almost onto. The 
property of being quasi-isometric is an equivalence relation on 
metric spaces.

\subsection{Subsurface projections and complexes of curves}
\label{surface complex background}

On any surface $S$, one may consider the 
\emph{complex of curves of $S$}, denoted $\CC(S)$. The complex of 
curves is a finite dimensional flag complex whose vertices correspond 
to nontrivial homotopy classes of non-peripheral, simple, closed curves and with 
edges between any pair of such curves which can be realized 
disjointly on $S$. In the cases where $r(S)\leq 1$ the definition 
must be modified slightly. When  $S$ is a one-holed torus or 4-holed
sphere, any pair of curves 
intersect, so edges are placed between 
any pair of curves which realize the minimal possible intersection on 
$S$ (1 for the torus, 2 for the sphere). 
With this modified definition, these curve complexes are the 
Farey graph. 
When $S$ is the 3-holed  sphere 
its curve complex is empty since $S$ supports no simple closed 
curves. Finally, the case when $S$ is an annulus 
will be important when $S$ is a subsurface of a larger surface $S'$. 
We define $\CC(S)$ by considering the annular cover
$\tilde{S'}$ of $S'$ in which $S$ lifts 
homeomorphically. Now $\tilde S'$ has a natural compactification to a
closed annulus, and
we let vertices be paths connecting the boundary components of this
annulus, up to homotopy rel endpoints. Edges are pairs of
paths with disjoint interiors. 
With this definition, one 
obtains a complex quasi-isometric to $\Z$. (See 
\cite{MasurMinsky:complex1} for further details.)

The following basic result on the curve complex was proved
by Masur--Minsky \cite{MasurMinsky:complex1}. (See also Bowditch
\cite{Bowditch:Curvecomplex}
for an alternate proof).
\begin{theorem}{} For any surface $S$, 
    the complex of curves is an infinite diameter $\delta$-hyperbolic 
    space (as long as it is non-empty).
\end{theorem}

Given a subsurface $Y\subset S$, one can define a \emph{subsurface 
projection} which is a map $\pi_{\CC(Y)}\co \CC(S)\to 2^{\CC(Y)}$. 
Suppose first $Y$ is not an annulus. Given any 
curve $\gamma\in\CC(S)$ intersecting $Y$ essentially, 
we define 
$\pi_{\CC(Y)}(\gamma)$ to be the collection of vertices in $\CC(Y)$ 
obtained by surgering the essential arcs of $\gamma\cap Y$ along 
$\boundary Y$ to obtain simple closed 
curves in $Y$. It is easy to show that 
$\pi_{\CC(Y)}(\gamma)$ is non-empty and has uniformly bounded diameter. 
If $Y$ is an annulus and $\gamma$ intersects it transversely
essentially, we may 
lift $\gamma$ to an arc crossing the annulus $\tilde{S'}$ and let this
be $\pi_{\CC(Y)}(\gamma)$. If $\gamma$ is a core curve of $Y$ or fails
to intersect it, we let $\pi_{\CC(Y)}(\gamma) = \emptyset$ (this holds
for general $Y$ too).

When measuring distance in the image subsurface, we usually write 
$d_{\CC(Y)}(\mu,\nu)$ as shorthand for 
$d_{\CC(Y)}(\pi_{\CC(Y)}(\mu),\pi_{\CC(Y)}(\nu))$.

\medskip 

\subsubsection*{Markings}

The curve complex can be used to produce a geometric model for the 
mapping class group as done in \cite{MasurMinsky:complex2}. This 
model is a graph called the \emph{marking complex}, 
$\MM(S)$, and is defined as follows. 

We define vertices $\mu\in\MM(S)$ to be pairs 
$(\base(\mu), \mbox{transversals})$ for which:
\begin{itemize}
    \item The set of \emph{base curves of $\mu$}, denoted $\base(\mu)$, is a 
    maximal simplex in $\CC(S)$.

    \item The \emph{transversals of $\mu$} consist of one curve for
    each component of $\base(\mu)$, intersecting it transversely.
\end{itemize}
Further, the markings are required to satisfy the following two 
properties. First, for each $\gamma\in\base(\mu)$, we require the 
transversal curve to $\gamma$, denoted $t$, be disjoint from the rest of the 
$\base(\mu)$. Second, given $\gamma$ and its transversal $t$, we 
require that $\gamma\cup t$ fill a non-annular surface $W$ satisfying 
$r(W)=1$ and for which $d_{\CC(W)}(\gamma, t)=1$.

The edges of $\MM(S)$ are of two types:
\begin{enumerate}
    \item \emph{Twist}: Replace a transversal curve by another 
    obtained by performing a Dehn twist along the associated base 
    curve.

    \item \emph{Flip}: Swap the roles of a base curve and its
    associated transversal curve.  
    (After doing this move, the additional disjointness requirement on
    the transversals may not be satisfied.  As shown in
    \cite{MasurMinsky:complex2}, one can surger the new transversal 
    to obtain one that does satisfy the disjointness requirement.
    The additional condition that the new and old transversals
    intersect minimally restricts the surgeries to a finite number, 
    and we obtain a finite set of possible flip moves for each marking.
    Each of these moves gives rise to an edge in the marking graph, and the
    naturality of the construction makes it invariant by the mapping
    class group.)
\end{enumerate}

It is not hard to verify that $\MM(S)$ is a locally finite graph on 
which the mapping class group acts cocompactly and properly 
discontinuously.  As observed by Masur--Minsky 
\cite{MasurMinsky:complex2}, this yields:

\begin{lemma}{} $\MM(S)$ is quasi-isometric to the mapping class 
    group of $S$.
\end{lemma}

The same definitions apply to essential subsurfaces of $S$. For an annulus $W$,
we let $\MM(W)$ just be $\CC(W)$.

Note that the above definition of marking makes no requirement 
that the surface $S$ be connected. In the case of a disconnected 
surface $\displaystyle W=\sqcup_{i=1}^{n} W_{i}$, 
it is easy to see that $\MM(W)=\prod_{i=1}^{n}\MM(W_{i})$.

\medskip

\subsubsection*{Projections and distance}
We now recall several ways in which subsurface projections arise in
the study of mapping class groups.

First, note that for any $\mu\in\MM(S)$ and any $Y\subseteq S$ 
the above projection maps extend to  $\pi_{\CC(Y)}\co\MM(S)\to2^{\CC(Y)}$. 
This map is simply the union over $\gamma\in\base(\mu)$ of the 
usual projections $\pi_{\CC(Y)}(\gamma)$, unless $Y$ is an annulus about an 
element of $\base(\mu)$. When $Y$ is an annulus about 
$\gamma\in\base(\mu)$, then we let $\pi_{\CC(Y)}(\mu)$ be the 
projection of $\gamma$'s transversal curve in $\mu$. As in the case of curve 
complex projections, we write $d_{\CC(Y)}(\mu,\nu)$ as shorthand for 
$d_{\CC(Y)}(\pi_{\CC(Y)}(\mu),\pi_{\CC(Y)}(\nu))$. 

\begin{remark}\label{shared base curve} An easy, but useful, fact is that 
    if a pair of markings $\mu,\nu\in\MM(S)$ share a base curve 
    $\gamma$ and $\gamma\cap Y\neq\emptyset$, then there is a 
    uniform bound on the diameter of 
    $\pi_{\CC(Y)}(\mu)\cup\pi_{\CC(Y)}(\nu)$.
\end{remark}

We say a pair of subsurfaces \emph{overlap} if 
they intersect, and neither is nested in the other.
The following is proven in \cite{Behrstock:asymptotic}:
\begin{theorem}{big implies small}
    Let $Y$ and $Z$ be a pair of 
    subsurfaces of $S$ which overlap.
    There exists a constant $M_{1}$ depending only on the topological 
    type of $S$, such that for any $\mu\in\MM(S)$:
    $$\min\left\{d_{\CC(Y)}(\boundary Z, \mu),
    d_{\CC(Z)}(\boundary Y, \mu)\right\} \le M_{1}.
$$
\end{theorem}

Another application of the projection maps 
is the following distance formula of 
Masur--Minsky \cite{MasurMinsky:complex2}:
\begin{theorem}{distance formula} 
    If $\mu,\nu\in \MM(S)$,  then there exists a constant $K(S)$, 
    depending only on the topological type of $S$, such that for each 
    $K>K(S)$ there exists $a\geq 1$ and $b\geq 0$ for which:
    $$
    d_{\MM(S)}(\mu,\nu) \approx_{a,b} \sum_{Y\subseteq S} 
    \Tsh K{d_{\CC(Y)}(\pi_{\CC(Y)}(\mu),\pi_{\CC(Y)}(\nu))}
    $$
\end{theorem}
Here we define the expression
$\Tsh K{N}$ to be $N$ if $N>K$ and $0$ otherwise -- hence $K$
functions as a ``threshold'' below which contributions are ignored.

\subsubsection*{Hierarchy paths}
In fact, the distance formula of Theorem \ref{distance formula} is 
a consequence of a construction in 
\cite{MasurMinsky:complex2} of a class of quasi-geodesics in $\MM(S)$
which we call {\em hierarchy paths}, and which have the following
properties. 

Any two points $\mu,\nu\in\MM(S)$ are connected by at least one
hierarchy path $\gamma$. Each hierarchy path is a quasi-geodesic, with
constants depending only on the topological type of $S$. 
The path 
$\gamma$ ``shadows'' a $\CC(S)$-geodesic
$\beta$ joining $\base(\mu)$ to $\base(\nu)$, in the following sense:
There is a monotonic map $v\co\gamma\to\beta$, such that
$v(\gamma_n)$ is a vertex in $\base(\gamma_n)$ for every $\gamma_n$ in $\gamma$.

(Note: the term ``hierarchy'' refers to a long combinatorial
construction which yields these paths, and whose details we will not
need to consider here).

Furthermore the following criterion constrains the makeup of these
paths. It asserts that subsurfaces of $S$ which ``separate'' $\mu$
from $\nu$ in a significant way must play a role in the hierarchy
paths from $\mu$ to $\nu$:

\begin{lemma}{large links}
There exists a constant $M_{2}=M_{2}(S)$ such that, if 
$W$ is an essential
subsurface of $S$ and $d_{\CC(W)}(\mu,\nu) > M_{2}$, then for any
hierarchy path $\gamma$ connecting $\mu$ to $\nu$, there exists
a marking $\gamma_n$ in $\gamma$ with $[\boundary W] \subset
\base(\gamma_n)$. Furthermore there exists a vertex $v$ in the geodesic
$\beta$ shadowed by $\gamma$ such that $W\subset S\setminus v$. 
\end{lemma}
This follows directly from Lemma 6.2 of \cite{MasurMinsky:complex2}.

\subsubsection*{Marking projections}
We have already defined two types of subsurface projections; we end 
by mentioning one more which we shall use frequently. Given a 
subsurface $Y\subset S$, we define a projection 
$$
   \pi_{\MM(Y)}\co\MM(S)\to \MM(Y)
$$
using the following procedure:
If $Y$ is an annulus $\MM(Y) = \CC(Y)$, we let $\pi_{\MM(Y)} =
\pi_{\CC(Y)}$. For nonannular $Y$: given a marking $\mu$ we
intersect its base curves with $Y$ and choose a curve $\alpha\in
\pi_Y(\mu)$. We repeat the construction with the subsurface
$Y\setminus\alpha$, continuing until we have found a maximal simplex
in $\CC(Y)$. This will be the base of $\pi_{\MM(Y)}(\mu)$. The
transversal curves of the marking are obtained by projecting $\mu$ to
each annular complex of a base curve, and then choosing a transversal
curve which minimizes distance in the annular complex to this
projection. (In case a base curve of $\mu$ already lies in $Y$, this
curve will be part of the base of the image, and the transversal curve
in $\mu$ will be used to determine the transversal for the image). 

This definition involved arbitrary choices, but it is shown in
\cite{Behrstock:asymptotic} that the set of all possible choices form
a uniformly bounded diameter subset of $\MM(Y)$. 
Moreover, it is shown there that: 

\begin{lemma}{piM properties}
$\pi_{\MM(Y)}$ is coarsely Lipschitz with uniform constants. 
\end{lemma}
Similarly to the case 
of curve complex projections, we write 
$d_{\MM(Y)}(\mu,\nu)$ as shorthand for 
$d_{\MM(Y)}(\pi_{\MM(Y)}(\mu),\pi_{\MM(Y)}(\nu))$.

\subsection{Asymptotic cones}
\label{cone defs}
The asymptotic cone of a metric space is roughly defined to be the limiting 
view of that space as seen from an arbitrarily large distance. 
This can be made precise using ultrafilters:

By a \emph{(non-principal) ultrafilter} we mean a finitely
additive probability measure $\omega$ defined on the power set of the
natural numbers and taking values only $0$ or $1$, 
and for which every finite set has zero measure. The existence of
non-principal ultrafilters depends in a fundamental way on the Axiom
of Choice.  

Given a
sequence of points $(x_{n})$ in  a topological space $X$,
we say $x\in X$ is its \emph{ultralimit}, or
$x=\ulim x_{n}$, if for every
neighborhood $U$ of $x$ the set 
$\{n:x_n \in U\}$ has $\omega$-measure equal to 1.
We note that ultralimits are unique when they exist, and that when $X$
is compact every sequence has an ultralimit. 

The ultralimit of a 
sequence of based metric spaces $(X_{n}, x_{n}, \dist_{n})$ is defined
as follows: 
Using the notation $\seq y = (y_n\in X_n)\in \Pi_{n\in\N}X_n$
to denote a sequence, define
$\dist(\seq y,\seq z) = \ulim(y_n,z_n)$, 
where the ultralimit is taken in the compact set $[0,\infty]$. 
We then let
$$
\ulim (X_{n},x_{n}, \dist_{n})\equiv\{\seq y: \dist(\seq y,\seq x) <
\infty\}/\sim,
$$
where we define $\seq y\sim \seq y'$ if 
$\dist(\seq y,\seq y')=0$. Clearly $\dist$ makes this quotient into a
metric space.

Given a sequence of positive 
constants  $s_{n}\to\infty$ and a sequence $(x_n)$ of basepoints in a
fixed metric space $(X,\dist)$, we may consider the rescaled space
$(X,x_n,\dist/s_n)$. The ultralimit of this sequence is 
called the \emph{asymptotic cone of $(X,\dist)$ relative to the 
ultrafilter $\omega$, scaling constants $s_{n}$, and basepoint 
$\seq x = (x_{n})$}:
$$\cone(X,(x_{n}),(s_{n}))=\ulim (X, x_{n}, \frac{\dist}{s_{n}}).$$
(For further details see \cite{DriesWilkie, Gromov:PolynomialGrowth}.)

\medskip

For the remainder of the paper, let us
fix a non-principal ultrafilter $\omega$, a sequence of scaling
constants $s_n\to \infty$, and a basepoint $\mu_0$
for $\MM(S)$.
We write $\AM=\AM(S)$ to denote an 
asymptotic cone of  $\MM(S)$ with respect to these 
choices.
Note that since $\MM$ is quasi-isometric to a word metric on $\MCG$, 
the space $\AM$ is homogeneous and thus the asymptotic cone is 
independent of the choice of basepoint. Further, since on a given 
group any two finitely generated word metrics 
are quasi-isometric, fixing an ultrafilter and scaling constants we 
have that 
different finitely generated word metrics on $\MCG$ have bi-Lipschitz 
homeomorphic asymptotic cones. 
Also, we note that in general 
the asymptotic cone of a geodesic space is a geodesic space. 
Thus, $\AM$ is a geodesic space, and in particular is locally path connected.

Any essential connected subsurface $W$ inherits a basepoint
$\pi_{\MM(W)}(\mu_0)$, canonical up to bounded error by
Lemma~\ref{piM properties}, 
and we can use this to define its asymptotic cone
$\AM(W)$. For a disconnected subsurface
$W=\sqcup_{i=1}^{k}W_{i}$ we have $\MM(W)=\Pi_{i=1}^k\MM(W_i)$ and we may
similarly construct $\AM(W)$ which can be identified
with $\Pi_{i=1}^{k}\AM(W_{i})$ (this follows  
from the general fact that the process of taking asymptotic 
cones commutes with finite products). 
Note that for an annulus $A$ we've defined $\MM(A) =\CC(A)$ which is
quasi-isometric to $\Z$, so $\AM(A)$ is $\R$.

It will be crucial to generalize this to {\em sequences} of
subsurfaces in $S$. Let us note first the general fact that
any sequence in a finite set $A$ is $\omega$-a.e.\ constant.
That is, given $(a_n\in A)$ there is a unique $a\in A$ such that 
$\omega(\{n:a_n = a\}) = 1$.
Hence for example if $\seq W =
(W_n)$ is a sequence of essential subsurfaces of $S$ then the
topological type of $W_n$ is $\omega$-a.e.\ constant and we call this
the topological type of $\seq W$. Similarly the topological type of
the pair $(S,W_n)$ is $\omega$-a.e.\ constant. 
We can  moreover interpret expressions like
$\seq U \subset \seq W$ for sequences $\seq U$ and $\seq W$ of
subsurfaces to mean $U_n\subset W_n$ for $\omega$-a.e.\ $n$, and so
on. We say that two sequences $(\alpha_n)$, $(\alpha'_n)$
are equivalent mod~$\omega$ if 
$\alpha_n = \alpha'_n$ for $\omega$-a.e.\ $n$, and note that
topological type, containment etc are invariant under this equivalence
relation. 
Throughout, we adopt the convention of using boldface to denote 
sequences. We will always consider such sequences mod~$\omega$, unless
they are sequences of markings $\seq \mu \in \AM$, in which case they
are considered modulo the weaker equivalence $\sim$ from the
definition of asymptotic cones.

If $\seq W=(W_n)$ is a sequence of subsurfaces, we let 
$\AM(\seq W)$ denote the ultralimit of 
$\MM(W_n)$ with metrics rescaled by $\frac{1}{s_{n}}$ 
and with basepoints $\pi_{\MM(W_{n})}(\mu_0)$.
Note that $\AM(\seq W)$ can be 
identified with $\AM(W)$, where $W$ is a surface 
homeomorphic to $W_{n}$ for $\omega$-a.e.~$n$.

\section{Product regions}
\label{Qregions}
In this section we will describe the geometry of the set of markings
containing a prescribed set of base curves. Equivalently, in the
mapping class group such a set corresponds to the coset of the
stabilizer of a simplex in the complex of curves. Not surprisingly,
these regions coarsely decompose as products.

Let $\Delta$ be a simplex in the complex of curves, i.e., a multicurve
 in $S$. We may partition $S$ into
subsurfaces isotopic to complementary components of $\Delta$, and
annuli whose cores are elements of $\Delta$. After throwing away
components homeomorphic to $S_{0,3}$ we obtain what we call the
 ``partition'' of $\Delta$, and denote $\sigma(\Delta)$.

Let $\QQ(\Delta)\subset \MM(S)$ denote the set of markings whose bases
contain $\Delta$. There is a natural (coarse) identification
\begin{equation}\label{Q product}
\QQ(\Delta) \approx \prod_{U\in \sigma(\Delta)} \MM(U)
\end{equation}
where if $U$ is an annulus we take $\MM(U)$ to mean the annulus
complex of $U$.  This identification is obtained simply by restriction
(or equivalently by subsurface projection)
for each non-annulus component, and by associating transversals with 
points in annulus complexes for the annular components. 

Theorem~\ref{distance formula}  yields the
following basic lemmas. When $A$ is a subsurface and $B$ is a 
collection of curves, we write $\tint{A}{B}$ to mean that $B$
cannot be deformed away from $A$.

\begin{lemma}{Q QI product}
The identification (\ref{Q product}) is a quasi isometry with uniform
constants. 
\end{lemma}

\begin{lemma}{distance from Q}
If $\mu\in \MM(S)$ then 
$$
d(\mu,\QQ(\Delta)) \approx \sum_{\tint{W}{\Delta}} 
\Tsh K{d_{\CC(W)}(\mu,\Delta)}
$$
\end{lemma}

\begin{proof}[Proof of Lemma \ref{Q QI product}]
If $\mu,\nu\in \QQ(\Delta)$, the distance formula in
Theorem~\ref{distance formula} gives 
$$
d(\mu,\nu) \approx \sum_{W} \Tsh K{d_{\CC(W)}(\mu,\nu)}
$$
where the constants in $\approx$ depend on the threshold $K$.
Now if $\tint{W}{\Delta}$, then Remark~\ref{shared base curve} 
implies that $\pi_W(\mu)$ and $\pi_W(\nu)$ are each a
bounded distance from $\pi_W(\Delta)$, 
and hence the $W$ term in the
sum is bounded by twice this. Raising $K$ above this constant means
that all such terms vanish and the sum is only over surfaces $W$
disjoint from $\Delta$, or annuli whose cores are components of
$\Delta$. But this is estimated by  the distance in
$\prod_{U\in\sigma(\Delta)}\MM(U)$, using Theorem \ref{distance
  formula} in each $U$ separately.
\end{proof}

\begin{proof}[Proof of Lemma \ref{distance from Q}]
Let $\mu\in\MM(S)$. For any $\nu\in \QQ(\Delta)$, we note that, if
$\tint{W}{\Delta}$, then 
$$|d_{\CC(W)}(\mu,\nu) - d_{\CC(W)}(\mu,\Delta)| \le c$$
for some constant $c$, by Remark~\ref{shared base curve}.
If $K_0$ is the minimal threshold that can be used in the distance
formula of Theorem \ref{distance formula}, 
let $K = K_0 + 2c$. 
We then see that for any $W$ contributing to  the sum
$$
\sum_{\tint{W}{\Delta}} \Tsh K{d_{\CC(W)}(\mu,\Delta)}
$$
we must have 
$$
d_{\CC(W)}(\mu,\nu) \ge  d_{\CC(W)}(\mu,\Delta) - c >  K_0 
$$
and, since our choice of $K$ yields $\half d_{\CC(W)}(\mu,\Delta)>c$, 
we furthermore have
$$
d_{\CC(W)}(\mu,\nu) \ge \half d_{\CC(W)}(\mu,\Delta).  
$$
It follows then that
\begin{align*}
\sum_{W} \Tsh {K_0}{d_{\CC(W)}(\mu,\nu)} & \ge
\sum_{\tint{W}{\Delta}} \Tsh {K_0}{d_{\CC(W)}(\mu,\nu)} \\
& \ge
\half \sum_{\tint{W}{\Delta}} \Tsh K{d_{\CC(W)}(\mu,\Delta)}.
\end{align*}
This gives one direction of the desired inequality.

To obtain the other direction, we fix $\mu\in\MM(S)$ and 
let $\nu\in \QQ(\Delta)$ be the marking whose restriction
to each $U\in\sigma(\Delta)$ is just $\pi_{\MM(U)}(\mu)$. 
With this choice, 
$$
d_{\CC(W)}(\mu,\nu) \le c
$$
for a uniform constant $c$ whenever $\notint{W}{\Delta}$, since
the intersections of $\mu$ and $\nu$ with $W$ are essentially the same.
Setting our
threshold  $K\ge K_0+2c$ again these terms all vanish, and
\begin{align*}
\sum_{W} \Tsh K{d_{\CC(W)}(\mu,\nu)} &=
\sum_{\tint{W}{\Delta}} \Tsh K{d_{\CC(W)}(\mu,\nu)} \\
& \le 2 \sum_{\tint{W}{\Delta}} \Tsh {K_0}{d_{\CC(W)}(\mu,\Delta)}
\end{align*}
where the last inequality is obtained using the same threshold trick
as above (we can assume it is the same value of $c$).  
\end{proof}

\subsubsection*{Product regions in the asymptotic cone}
Consider a  sequence $\seq\Delta=\{\Delta_n\}$ 
such that $\ulim \frac{1}{s_n}d(\mu_0, \QQ(\Delta_n)) < \infty$. 
We can take the
ultralimit of $\QQ(\Delta_n)$, with  metrics rescaled by $1/s_n$,
obtaining a subset of $\AM(S)$ 
which we denote $\AQ(\seq\Delta)$.
Lemma \ref{Q QI product} and the fact that ultralimits commute 
with finite products implies that there is a bi-Lipschitz identification 
\begin{equation}
\label{Q product structure}
\AQ(\seq\Delta) \homeo 
\prod_{\subseq U \in \sigma(\subseq\Delta)} \AM(\seq
U). 
\end{equation}
Here $\sigma(\seq\Delta)$ is defined as follows: As in Section
\ref{cone defs},
the topological type of $\sigma(\Delta_n)$ is $\omega$-a.e.\ constant,
and so there is a set $J\subset \N$ with $\omega(J)=1$, 
a partition $\sigma' = \{U^1,\ldots,U^k\}$ of $S$, and a sequence of
homeomorphisms $f_n\co S\to S$ taking $\sigma'$ to $\sigma(\Delta_n)$ for
each $n\in J$. We then let $\sigma(\seq\Delta) = \{\seq
U^1,\ldots,\seq U^k\}$ where $\seq U^i = (f_n(U^i))$ for $n\in J$ (it
doesn't matter, mod~$\omega$,  how we define it for $n\notin J$). Any
non-uniqueness of $f_n$, up to isotopy, corresponds to a symmetry of
$\sigma'$, and hence to a permutation of the indices of elements of 
$\sigma(\seq\Delta)$.

Moreover, Lemma \ref{distance from Q} implies that
distance to $\AQ(\seq\Delta)$ can be estimated, up to bounded ratio,
by:
\begin{equation}\label{rho define}
\rho(\seq\mu,\Delta) \equiv
\lim_\omega \frac{1}{s_n} \sum_{\tint{W}{\Delta_n}}
  \Tsh K{d_{\CC(W)}(\mu_n,\Delta_n)}.
\end{equation}

\section{Separating product regions and locally constant maps}
\label{Pregions}
In this section we will define the family of product regions equipped with
locally constant maps (denoted as $\PP$ in the outline in the
introduction). Each region will be determined by a sequence $\seq W
=(W_n)$ of connected subsurfaces of $S$, and a choice $\seq x = (x_n)$
of basepoint in $\AM(\seq W)$.  Theorem \ref{global projection to F},
which defines the projection map associated to each region and
establishes its properties, 
is the main result of this section.

\subsection{Sublinear growth sets}
In Behrstock \cite{behrstock:thesis}, a family of subsets of $\AM(S)$ is
introduced, and defined as follows: 
for $\seq x\in\AM(S)$, let 
$$
F(\seq x) = \left\{\seq y: 
\lim_\omega \frac{1}{s_n} \sup_{U\subsetneq S} d_{\MM(U)}(x_n,y_n) = 0
\right\}.
$$
That is, the distance between $x_n$ and $y_n$, projected to the
marking graph of any proper subsurface, is vanishingly small compared to
their distance in $\MM(S)$.  We note that, because the subsurface
projections are uniformly Lipschitz, this condition is well-defined,
i.e., does not depend on the choice of $y_n$ representing $\seq y$. 

Behrstock proved that $F(\seq x)$ is an $\R$-tree, and more strongly that
for any two points in $F(\seq x)$ there is a unique embedded arc in $\AM(S)$
connecting them. We can generalize this construction slightly as
follows:

First, for a sequence $\seq U = (U_n)$ of connected subsurfaces and
$\seq x,\seq y \in \AM(S)$ we have 
$$
d_{\AM(\subseq U)}(\seq x, \seq y) = \lim_\omega \frac{1}{s_n} 
d_{\MM(U_n)}(x_n,y_n). 
$$
Now if $\seq W=(W_n)$ is a sequence of connected subsurfaces (considered
mod $\omega$) and $\seq x\in \AM(\seq W)$,  
we define $F_{\subseq W,\subseq x}\subset \AM(\seq W)$ to be:
$$
F_{\subseq W,\subseq x} = \{ \seq y\in\AM(\seq W): d_{\AM(\subseq U)}(\seq x,\seq y)=0\ 
\  \text{for all   $\seq U\subsetneq \seq W$} \}.
$$
If $W_n \equiv S$, this is equivalent to the definition of $F(\seq x)$
above. 
Note also that if
$\seq W=\collar(\balpha)$ then $F_{\subseq{W,x}}$ is just the 
asymptotic cone of the annulus complex of $\seq W$, 
which is a copy of $\R$.

Let us restate and discuss Behrstock's theorem from 
\cite{Behrstock:asymptotic}:
\begin{theorem}{F projection}
Let $\seq W = (W_n)$ be a sequence of connected subsurfaces of 
$S$, and
$\seq x \in \AM(\seq W)$. 
Any two points $\seq y,\seq z\in F_{\subseq W,\subseq x}$ are connected by a
unique embedded path in 
$\AM(\seq W)$, and this path lies in $F_{\subseq W,\subseq x}$. 
\end{theorem}

In particular, it follows that $F_{\subseq{W,x}}$ is an $\R$-tree. 

Here is a brief outline of the proof: 
The annular case is trivial because 
$F_{\subseq W, \subseq x} = \AM(\seq W) \homeo \R$.
Hence, we assume $W_n$ are not annuli for $\omega$-a.e.~$n$.
In each $W_n$, connect $y_n$ to
$z_n$ with a hierarchy path $\gamma_n$ (see \S\ref{surface complex
  background}). Since $\gamma_n$ are uniform quasi-geodesics, after
rescaling their ultralimit gives a path $\seq \gamma$ in $\AM(\seq W)$.
Using the tools of \cite{MasurMinsky:complex2} 
together with the assumption that $\seq y,\seq z \in F_{\subseq{W,x}}$, 
one can show that $\seq \gamma$
lies in $F_{\subseq{W,x}}$. 

Let $\beta_n$ be a $\CC(W_n)$-geodesic shadowed by $\gamma_n$.
One can see that the length $|\beta_n| \to_\omega \infty$ as follows: Suppose
instead that $|\beta_n|< L$ for $\omega$-a.e.~$n$. 
Choose the threshold in the distance formula large enough
so that the non-zero terms in 
$$
\sum_{V\subset W_n} \Tsh K{d_{\CC(V)}(y_n,z_n)}
$$
are proper subsurfaces in $W_n$ which play the role in $\gamma_n$
determined by Lemma \ref{large links} --- that is, each one is disjoint
from some $v\in\beta_n$. But since $\beta_n$ has at most $L$
vertices,  There must be one, $v_n$, which is disjoint from enough
surfaces to contribute at least $1/L$ times the sum. But this
means, using the distance formula within $Y_n = S\setminus v_n$, that
$d_{\AM(\subseq Y)}(\seq y,\seq z) > 0$, which contradicts the
assumption that $\seq y,\seq z \in F_{\subseq W,\subseq x}$.

Consider the map $p_n\co\MM(W_n) \to \beta_n$ which takes a
marking $\mu$ to a vertex $v\in \beta_n$ of minimal 
$\CC(W_n)$-distance to the base of $\mu$. 
We promote $p_n$ to a map 
$q_n\co \MM(W_{n})\to \gamma_n$
by letting $q_n(\mu)$ be a marking of $\gamma_n$ which shadows
$v=p_n(\mu)$. 

The ultralimit of $q_n$ yields a map $\seq q\co\AM(\seq W) \to
\seq\gamma\subset F_{\subseq W,\subseq x}$. Furthermore one can show
using hyperbolicity of $\CC(W_n)$ (Masur--Minsky
\cite{MasurMinsky:complex1}) and properties of the subsurface
projection maps that $q_n$ has coarse contraction properties that, in
the limit, imply that $\seq q$ is locally constant in the complement
of ${\seq \gamma}$.  It then easily
follows that $\seq y$ and $\seq z$ cannot be connected in the
complement of any point of $\seq \gamma$, and hence any path
between them must contain $\seq \gamma$, and any embedded path
must equal $\seq \gamma$.

\medskip

\subsection{Definition of $P_{\subseq W,\subseq x}$}

Given $\seq W $ and $\seq x$ as above, 
our separating product regions, denoted  $P_{\subseq W,\subseq x}$, will be
subsets of $\AQ(\boundary \seq W)$ defined as follows: 

In the product structure (\ref{Q product structure}) for
$\AQ(\boundary \seq W)$, $\seq W$ is a member of $\sigma(\boundary
\seq W)$, and hence $\AM(\seq W)$ 
appears as a factor. We let $P_{\subseq W,\subseq x}$ be the subset of
$\AQ(\boundary \seq W)$ consisting of points whose coordinate in
the $\AM(\seq W)$ factor lies in $F_{\subseq W,\subseq x}$.

Since the identification of $\AQ(\boundary \seq W)$ with the product structure
is made using the subsurface projections, we have this
characterization: 

\begin{lemma}{PW def}
$P_{\subseq W,\subseq x}$ is the set of points $\seq y\in\AM(S)$
such that:
\begin{enumerate}
\item $\pi_{\AM\subseq W}(\seq y) \in F_{\subseq W,\subseq x}$, and
\item $\rho(\seq y, \boundary\seq W) = 0.$
\end{enumerate}
\end{lemma}

Here $\rho(\seq y,\boundary\seq W)$ is an estimate for
the distance of $\seq y$ from $\AQ(\boundary\seq W)$, as defined
in (\ref{rho define}).
Also, the ultralimit of the
rescaled marking projection maps $\MM(S) \to \MM(W_n)$ is denoted by:
$$\pi_{\AM\subseq W}\co\AM(S) \to \AM(\seq W).$$ 
\noindent

Define $W_n^c$ to be the union of the components of 
$\sigma(\boundary W_n)$ not
equal to $W_n$ (so $W_n^c$ includes annuli around $\boundary W_n$,
unless $W_n$ itself is an annulus). Let $\seq W^c = (W_n^c)$. 
Then $\AM(\seq W^c)$ is the asymptotic cone of $(\MM(W_n^c))$, 
and can be identified with the product of the remaining factors in
$\AQ(\boundary \seq W)$: 
$$
\AM(\seq W^c) \equiv
\prod_{\substack{\subseq U\in \sigma(\boundary\subseq W) \\ U\ne \subseq W}}
\AM(\seq U)
$$
We can summarize this in the following: 

\begin{lemma}{product structure}
There exists a  bi-Lipschitz identification of
$P_{\subseq W,\subseq x}$ with
$$
F_{\subseq W,\subseq x} \times \AM(\seq W^c).
$$
\end{lemma}

\subsection{Projection maps}

The following projection theorem is a small improvement on 
Theorem \ref{F projection} from Behrstock \cite{Behrstock:asymptotic}.

\begin{theorem}{local projection to F}
Given $\seq x\in\AM(\seq W)$, 
there is a continuous map 
$$\wp=\wp_{\subseq W,\subseq x}\co \AM(\seq W) \to F_{\subseq W,\subseq x}$$
 with
these properties: 
\begin{enumerate}
\item $\wp$ is the identity on $F_{\subseq W,\subseq x}$
\item $\wp$ is locally constant in $\AM(\seq W) \setminus F_{\subseq
  W,\subseq x}$. 
\end{enumerate}
\end{theorem}

Note that in the proof of Theorem \ref{F projection} a projection to
individual paths was shown to have locally-constant properties. In
this theorem we construct a projection from $\AM(\seq W)$ onto $F_{\subseq{W,x}}$. 

\begin{proof}
For any $\seq y\in \AM(\seq W)$ let $\alpha$ be a path
connecting $\seq y$ to any point in $F_{\subseq W, \subseq x}$. Let
$\alpha_1$ be the first 
point in $\alpha$ that is in $F_{\subseq W, \subseq x}$. We claim that $\alpha_1$
depends only on $\seq y$. For otherwise let $\beta$ be another path with
$\beta_1\ne \alpha_1$. Then segments of $\alpha$ and $\beta$ form a
path connecting two points of $F_{\subseq W, \subseq x}$ outside of $F_{\subseq
  W, \subseq x}$ --- this contradicts
Theorem  \ref{F projection}.  

We can then define $\wp(\seq y) \equiv \alpha_1$. This is 
locally constant at $\seq y\notin F_{\subseq W, \subseq x}$ 
because for a sufficiently small
neighborhood $U$ of $\seq y$, every $\seq z\in U$ can be 
connected to $F_{\subseq W, \subseq x}$ by a path going 
first through $\seq y$ (since $\AM(\seq W)$ is 
locally path-connected). 

Continuity of  $\wp$ at points of  
$F_{\subseq W,\subseq x}$ follows immediately from the 
definition of $\wp$ and the fact that $\AM(\seq W)$ is a locally path 
connected geodesic space.
\end{proof}

We can now construct our global projection map for $F_{\subseq W, \subseq x}$: 

\begin{theorem}{global projection to F}
Given $\seq x\in\AM(\seq W)$, there is a continuous map
$$
\Phi=\Phi_{\subseq W, \subseq x} \co \AM(S) \to F_{\subseq W, \subseq x}
$$
with these properties:  
\begin{enumerate}
\item 
$\Phi$ restricted to $P_{\subseq W, \subseq x}$ 
is projection to the first factor in the product structure 
$P_{\subseq W, \subseq x} \homeo
  F_{\subseq W, \subseq x}\times \AM(\seq W^c)$.
\item
$\Phi$ is locally constant in the complement of $P_{\subseq W, \subseq x}$. 
\end{enumerate}
\end{theorem}

\begin{proof}
We define the map simply by 
$$
\Phi_{\subseq W, \subseq x}  = \wp_{\subseq W, \subseq x} \circ 
\pi_{\AM\subseq W}. 
$$
Property (1) follows from the definition, and from the way that the
identification of $P_{\subseq W, \subseq x}$ with the product in Lemma \ref{product
  structure} is constructed via subsurface projections.

We divide the proof of property (2) into two cases: 

{\bf Case 1:} $\pi_{\AM\subseq W}(\seq y)\notin F_{\subseq W,\subseq x}$. 

In this
case the desired 
fact follows immediately from the locally-constant property of $\wp$
shown in Theorem \ref{local projection to F}, and the continuity of
$\pi_{\AM\subseq W}$.

{\bf Case 2:} $\pi_{\AM\subseq W}(\seq y)\in F_{\subseq W,\subseq x}$.

Since $\seq y\notin P_{\subseq W,\subseq x}$ and 
$\pi_{\AM\subseq W}(\seq y)\in F_{\subseq W,\subseq x}$, 
Lemma~\ref{PW def} implies that $\rho(\seq y,\boundary \seq W)>0.$

Let $\seq z\in \AM(S)$, with $\Phi(\seq z) \ne \Phi(\seq y)$. We will
derive a lower 
bound for $d(\seq y,\seq z)$, and this will prove the theorem. 

Let $\seq z' = \pi_{\AM\subseq W}(\seq z)$ and $\seq y'=\pi_{\AM\subseq
  W}(\seq y)$. Since Case~1 has already been handled,  
  we may assume $\seq y'\in F_{\subseq W, \subseq x}$, so $\seq y' =
  \wp(\seq y') = \Phi(\seq y)$. 
 As in Theorem \ref{local projection to F},  any path from $\seq z'$ to $\seq y'$ must pass 
  through $\wp(\seq z')$ first.  Note that $\wp(\seq z') = \Phi(\seq  z) \neq \seq y'$.
Now let $\gamma_n$ be hierarchy paths in $\MM(W_n)$ connecting $z'_n$ to
  $y'_n$. Since $\gamma_n$ are quasigeodesics, their ultralimit after
  rescaling gives
  rise to a path in $\AM(\seq W)$ connecting $\seq z'$ to $\seq y'$ and
  hence there must exist
  $\delta_n\in \gamma_n$ such that $(\delta_n)$
  represents $\wp(\seq z')$. As remarked in the outline of the proof
  of Theorem \ref{F projection}, $d_{\CC(W_n)}(\delta_n,y'_n) \to_\omega \infty$ since
$\wp(\seq z')$ and $\seq y'$ are distinct points in $F_{\subseq
  W,\subseq x}$. Now since $\gamma_n$ monotonically shadows a
  $\CC(W_n)$ geodesic from $z'_n$ to $y'_n$, we conclude that
$$
d_{\CC(W_n)}(y'_n,z'_n) \to_\omega \infty.
$$
Since $\pi_{\CC(W_n)}\circ\pi_{\MM(W_n)}$ and $\pi_{\CC(W_n)}$
differ by a bounded constant (immediate from the
definitions), we conclude that
$$
d_{\CC(W_n)}(y_n,z_n) \to_\omega \infty. 
$$

Now by the definition of $\rho(\seq y,\boundary \seq W)$, we know that
\begin{equation}\label{d y W sum}
\frac{1}{s_n}\sum_{\tint{U}{\boundary W_n}}
\Tsh K{d_{\CC(U)}(y_n,\boundary W_n)}
\to_\omega  c>0.
\end{equation}
Let $U$ be a subsurface participating in this sum for some $n$, 
so that we have $d_{\CC(U)}(y_n,\boundary W_n) > K $. We want to show 
that 
\begin{equation}\label{U between}
d_{\CC(U)}(y_n,z_n) \ge d_{\CC(U)}(y_n,\boundary W_n) - K'
\end{equation}
for some $K'$.

We assume that $K$ is larger than the constant $M_{1}$ from  
Theorem~\ref{big implies small}, recall that this theorem states that 
\begin{equation}\label{J theorem}
\min \{ d_{\CC(V)}(\mu,\boundary V') , 
d_{\CC(V')}(\mu,\boundary V)\} \le M_{1} 
\end{equation}
for any  marking $\mu$ and
subsurfaces $V,V'$ with $\boundary V \pitchfork \boundary V' \ne
\emptyset$.

Since $U$ meets $\boundary W_n$, we have either $\tint{\boundary
  U}{W_n}$, in which case the subsurfaces $W_{n}$ and $U$ overlap, 
  or $W_n \subsetneq U$. 

Suppose first that $\tint{\boundary
  U}{W_n}$. 
Now we have $d_{\CC(U)}(y_n, \boundary W_n) > K>M_{1}$, since 
$W_{n}$ and $U$ overlap (\ref{J theorem}) implies 
$$
d_{\CC(W_n)}(y_n,\boundary U) \le M_{1}. 
$$
Now by the triangle inequality 
$$
d_{\CC(W_n)}(\boundary U, z_n) \ge d_{\CC(W_n)}(y_n,z_n) - M_{1} - D
$$
(where $D$ is a bound for $\diam_{\CC(W_n)}(\mu)$ of any marking, as 
given by Remark~\ref{shared base curve}). 
Since $d_{\CC(W_n)}(y_n,z_n) \to_\omega \infty$, we may assume that 
this gives
$$
d_{\CC(W_n)}(\boundary U,z_n) > M_{1}.
$$
Now again by (\ref{J theorem}) we have
$$
d_{\CC(U)}(\boundary W_n, z_n) \le M_{1}
$$
and again by the triangle inequality
$$
d_{\CC(U)}(y_n,z_n) \ge d_{\CC(U)}(y_n,\boundary W_n) - M_{1} -D 
$$
which establishes (\ref{U between}) when $\tint{\boundary
  U}{W_n}$.

Next, let us establish (\ref{U between}) when $W_n \subsetneq U$. 
Since $d_{\CC(W_n)}(y_n,z_n)\to_\omega \infty$, we may assume that
this distance is larger than the constant $M_{2}$ in Lemma~\ref{large links}. 
Let $\gamma_n$ be a hierarchy path in $\MM(U)$ connecting $\pi_{\MM
  U}(y_n)$ to $\pi_{\MM U}(z_n)$, and let $\beta_n$ be the
$\CC(U)$-geodesic from $\pi_{\CC(U)}(y_n)$ to $\pi_{\CC(U)}(z_n)$ that
$\gamma_n$ shadows.  Lemma~\ref{large links} implies that
$\boundary W_n$ appears in the base of at
least one marking in $\gamma_n$, and hence 
$[\boundary W_n]$ is $\CC(U)$-distance at most one from  
 a vertex of $\beta_n$. This means that the length of $\beta_n$
is at least $d_{\CC(U)}(\boundary W_n,y_n)-2$, in 
particular:
$$
d_{\CC(U)}(z_n,y_n) \ge d_{\CC(U)}(y_n,\boundary W_n) - 2.
$$
Thus, we have established (\ref{U between}) with $K'=\max\{M_{1}+D, 2\}$.

Now applying this to all the terms in the sum of (\ref{d y W sum}), we
would like to obtain a lower bound (for $\omega$-a.e.~$n$) 
\begin{equation}\label{sum lower bound}
\frac{1}{s_n}\sum_{\tint{U}{\boundary W_n}}
\Tsh K{d_{\CC(U)}(y_n,z_n)}
> c' > 0
\end{equation}
To do this we apply the same threshold trick we used in the proof of
Lemma \ref{distance from Q}. Since 
Theorem~\ref{distance formula} applies to any
sufficiently large threshold, we may choose $K'' = 2K'+K$ to 
replace the
threshold $K$ in the sum in (\ref{d y W sum}), and obtain
\begin{equation}\label{new threshold}
\frac{1}{s_n}\sum_{\tint{U}{\boundary W_n}}
\Tsh {K''}{d_{\CC(U)}(y_n,\boundary W_n)}
\to_\omega c' > 0.
\end{equation}
Now for a given $n$ if $U$ contributes to this sum then by
(\ref{U between}), we have
$d_{\CC(U)}(y_n,z_n) \ge K''-K' > K$, and moreover
$$
d_{\CC(U)}(y_n,z_n) \ge d_{\CC(U)}(y_n,\boundary W_n) - K' \ge \half
d_{\CC(U)}(y_n,\boundary W_n). 
$$
This implies that
$$
\sum_{\tint{U}{\boundary W_n}} 
\Tsh{K}{d_{\CC(U)}(y_n,z_n)} \ge \half
\sum_{\tint{U}{\boundary W_n}} 
\Tsh{K''}{d_{\CC(U)}(y_n,\boundary W_n)}.
$$
In other words, again 
using the distance formula, this gives us a lower bound of the form
$$
d_{\AM(S)}(\seq y,\seq z) > c'' > 0. 
$$
The conclusion is that if $d(\seq y,\seq z) < c''$ then $\Phi(\seq
y)=\Phi(\seq z)$, which is what we wanted. 
\end{proof}

\subsection{Separators}

In \cite{Behrstock:asymptotic}, it was shown that mapping class 
groups have global cut-points in their asymptotic cones, cf.\ 
Theorem~\ref{F projection}. Since mapping class groups are not 
$\delta$-hyperbolic, except in a few low complexity cases, it clearly 
can not hold that arbitrary pairs of points in the asymptotic cone  
are separated by a point. Instead we identify here a larger class of
subsets which do separate points: 

\begin{theorem}{separate points}
There is a family $\LL$ of closed subsets of $\AM(S)$ such that
any two points in $\AM(S)$ are separated by some $L\in\LL$. 
Moreover each $L\in \LL$ is isometric to
$\AM(Z)$, where $Z$ is some proper essential (not necessarily
connected) subsurface of $S$, with $r(Z) < r(S)$. 
\end{theorem}

We will see as 
part of an inductive argument
in the next section that
these separators $L$ all have  
(locally compact) dimension at most $r(S)-1$; this bound is sharp since $\AM$ 
contains $r(S)$-dimensional bi-Lipschitz flats which, of course, can 
not be separated by any subset of dimension less than $r(S)-1$.

\begin{proof} Fix $\seq x\ne \seq y\in \AM(S)$. We claim that there exists
a subsurface sequence $\seq W=(W_n)$ such that:
\begin{enumerate}
\item $d_{\AM(\subseq W)}(\seq x,\seq y) > 0$, and
\item For any $\seq Y=(Y_n)$ with $\seq Y\subsetneq \seq W$, 
$d_{\AM(\subseq Y)}(\seq x,\seq y) = 0$.
\end{enumerate}
Indeed, $\seq W=(S)$ satisfies the first condition. If it fails the second,
we may choose $\seq{W'}\subsetneq \seq W$ with 
$d_{\AM(\subseq{W'})}(\seq x,\seq y) > 0$, and 
continue. This terminates since the complexity of the subsurface
sequence decreases.

Let $\seq x' = \pi_{\AM\subseq W}(\seq x)$ 
and $\seq y' = \pi_{\AM\subseq W}(\seq y)$. 
The choice of $\seq W$ implies that $\seq x' \ne \seq y'$ and that 
$\seq y'\in F_{\subseq W, \seq x'}$. (Note that the second condition 
implies  $F_{\subseq W, \seq x'} = F_{\subseq W, \seq y'}$.)
Let $\seq z$ be a point in $F_{\subseq W,\seq x'}$ in the interior of 
the path from $x'$
to $\seq y'$. Since $F_{\subseq W,\seq x'}$ is an $\R$-tree (by 
Theorem~\ref{F projection}), $\seq z$ separates $\seq x'$ from 
$\seq y'$ in $F_{\subseq W,\seq x'}$.

Let $L$ be the subset of $P_{\subseq W,\seq x'}$ identified with
$\{\seq z\}\times\AM(\seq W^c)$ by Lemma \ref{product structure}.
Certainly $L$ separates $P_{\subseq W,\seq x'}$. We claim 
$L$ also separates
$\AM(S)$, with $\seq x$ and $\seq y$ on different sides. This follows
immediately from Theorem \ref{global projection to F}:

Recall the map
$\Phi=\Phi_{\subseq W,\seq x'}\co\AM(S) \to F_{\subseq W, \seq x'}$,  
and that $\seq x'=\Phi(\seq x)$ and $\seq y'=\Phi(\seq y)$. 
Divide $F_{\subseq W,\seq x'}\setminus\{\seq z\}$ into two disjoint  
open sets
$E_{\seq x}$ and $E_{\seq y}$ containing $\seq x'$ and $\seq y'$, 
respectively. 
$\Phi^{-1}(E_{\seq x})$ and $\Phi^{-1}(E_{\seq y})$ are open 
sets containing $\seq x$ and
$\seq y$ respectively. The remainder $\Phi^{-1}(\{\seq z\})$ 
consists of $L$
union an open set $V$, by the locally constant property.
Hence we have divided $\AM(S)\setminus L$ into three disjoint open
sets two of which contain $\seq x$ and $\seq y$ respectively. 
This proves $L$ separates $\seq x$ and $\seq y$.

The construction exhibits $L$ as an asymptotic cone $\AM(\seq W^c)$, 
from which it follows that $L$ is closed (cf. \cite{DriesWilkie}).
Since the topological type of $\seq W^c$ is $\omega$-a.e.~constant,
this is isometric to $\AM(W^c)$ for some fixed surface $W^c$. 
\end{proof}

\section{The dimension theorem}
\label{dimension}

In this section we will apply the separation Theorem \ref{separate
  points} to prove the main theorem on dimension in $\AM(S)$.
We begin with some terminology: 

Historically, topologists have studied three different versions of
dimension: small inductive dimension, $\ind$,
large inductive dimension, $\Ind$, and covering dimension, $\Dim$ 
(the covering dimension is also called the topological dimension).
Dimension theory grew out of the development of these various
definitions and studies the interplay and applications of the various
versions of dimension \cite{Engelking:DimBook}.  For a topological
space $X$, let $\indhat(X)$ denote the supremum of $\ind(X')$ over all
locally-compact subsets $X'\subset X$, and similarly define $\Indhat$
and $\dimhat$.
 Restating our main theorem, we have:

\begin{theorem}{inductive dimension}
$\indhat(\AM(S))=\Indhat(\AM(S))=\dimhat(\AM(S)) =  r(S)$.
\end{theorem}

The Rank Conjecture follows immediately as a corollary, since
$\R^n$ is locally compact and $\ind(\R^n) =n$. 

\subsection{Separation and dimension}

We will work with inductive dimension, which we define
below. Equivalence of the different dimensions in our setting is
provided by 

\begin{lemma}{dimensions equal} For a metric space $X$, 
  $\dimhat(X) = \indhat(X) = \Indhat(X)$. 
\end{lemma}

\begin{proof} 
This is essentially an appeal to the literature. 
First note the following standard topological facts:
	    \begin{enumerate}
		\item every metric space is paracompact;
		\item a locally compact space is paracompact if and only if it
		is strongly paracompact 
		\cite[Page 329]{engelking:TopBook}.
	    \end{enumerate}
Engelking shows
\cite[pg 220]{Engelking:DimBook} that
if $Y$ is a strongly paracompact metrizable
space, then $\ind(Y)=\Ind(Y)=\Dim(Y)$.
Thus, if $X'\subset X$ is a locally compact subset,
then $\ind(X')=\Ind(X')=\Dim(X')$.  Taking the supremum over
locally compact subsets finishes the proof. 
\end{proof}

To prove Theorem~\ref{inductive dimension} we provide a lemma 
reducing this result to Theorem~\ref{separate points}. First we recall 
the definition of the \emph{small inductive dimension}:
$\ind(\emptyset) = -1$ and for any $X$, $\ind(X) = n$ if $n$
is the smallest number such that for all $x\in X$ and 
neighborhood $V$
of $x$, there exists a neighborhood $x\in U \subset V$ such that
$\ind(\boundary U) \le n-1$. Here $\boundary U$ is the topological 
frontier of $U$ in $Y$.  (See \cite{Engelking:DimBook} for further
details.)

\begin{lemma}{separators implies dimension}
    If $X$ is a metric space for which every pair of points can be 
    separated by a closed subset $L\subset X$ with 
    $\indhat(L) \le D-1$, then $\indhat(X)=\Indhat(X)=\dimhat(X)\le D$.
\end{lemma}

\begin{proof}
By Lemma \ref{dimensions equal}, we may henceforth restrict our
attention to the small inductive dimension.

Let $X'$ be a locally compact subset of $X$.
Fixing $x\in X'$, consider any $\ep$-ball $B$ about $x$ in the 
induced metric on $X'$, where $\ep$ is assumed to be sufficiently 
small so that 
local compactness of $X'$ implies $\boundary B$ is compact. For any
$y\in \boundary B$, let $L$ be a closed separator of 
$x$ and $y$, with
$\indhat(L)\le D-1$, as provided by hypothesis. 
Since $X'$ is locally compact, 
$L'=X'\intersect L$ has $\ind(L')\le D-1$.
The separation property means that $X'\setminus L'$ is the 
union of a pair of disjoint open subsets of $X'$,  
$W_y$ and $V_y$,  
such that $x\in W_y$ and $y\in V_y$. 
Since $\boundary B$ is compact, we may extract a
finite subcover of the covering $\{V_y\}$ of 
$\boundary B$, which we relabel $V_1,\ldots,V_n$, with corresponding 
separators $L_1,\ldots,L_n$ and complementary $W_1,\ldots,W_n$. 
Then $\union L'_i$ separates $x$ from $\boundary
B$.  More precisely, let
${\mathcal W}=\intersect W_i$ and ${\mathcal V}=\union
V_i$. (In case $\boundary B=\emptyset$, let ${\mathcal W} = X'$ and
${\mathcal V}=\emptyset$.)
These are disjoint open sets with $x\in {\mathcal W}$,
$\boundary B \subset {\mathcal V}$, and
$\boundary \mathcal W \subset \union L'_i$. 

Now let $U = \mathcal W \intersect B$. This is an open set,
contained in $B$, whose boundary is contained in $\union L'_i$
(since it cannot meet $\boundary B$ which lies in $\mathcal V$). 
Since $\ind$ is preserved by finite unions and monotonic with respect
to inclusion, we have
$\ind(\boundary U) \le D-1$, which is what we wanted to prove. 
\end{proof}

\subsection{Proof of the dimension theorem}
We can now complete the proof of 
Theorem~\ref{inductive dimension}, by inducting on $r(S)$.

Note that the lower bound $\indhat(\AM(S)) \ge r(S)$ is immediate
since maximal abelian subgroups give quasi-isometrically 
embedded $r(S)$-flats \cite{FarbLubotzkyMinsky}. We now
prove the upper bound. 

When $r(S)=1$, $S$ is
$S_{1,1}$, $S_{0,4}$ or $S_{0,2}$.  The asymptotic cones for the first two are
the asymptotic cone for $SL(2,\Z)$ which is known to be an $\R$-tree. 
In the third case we really have in mind the annulus complex of an
essential annulus, for which the asymptotic cone is just $\R$.
Since $\indhat = 1$ is well known for $\R$-trees, the theorem holds in
this case. 

Theorem~\ref{separate points} provides for each $x,y\in\AM(S)$
a separator, $L$, which is homeomorphic to $\AM(W^c)$, where $W$ is
an essential subsurface of $S$. Since $r$ is additive over disjoint
unions and $r(W)\ge 1$, we have $r(W^c) \le r(S)-1$. Thus by induction
$\indhat(L) \le r(S)-1$. (We can apply the inductive hypothesis to each
component of $W^c$, and use subadditivity of $\ind$ over finite products, 
see \cite{Engelking:DimBook}, and additivity of $r$ over disjoint unions.)

Thus we have satisfied the hypotheses of 
Lemma~\ref{separators implies dimension} for $\AM(S)$, and 
Theorem~\ref{inductive dimension} follows.

\section{Teichm\"{u}ller space}
\label{pants}
In this section we deduce analogues of the results 
in the earlier sections for Teichm\"{u}ller space with the 
Weil-Petersson metric. As shown in  Brock \cite{Brock:wp}, there is a 
combinatorial model for the Weil-Petersson metric on 
Teichm\"{u}ller space provided by the \emph{pants graph}. 
The combinatorial analysis as carried out above for the mapping class 
group can be done similarly in the pants graph, 
(cf. \cite[Section 8]{MasurMinsky:complex2}). Using Brock's result, we 
deduce the results below about Teichm\"{u}ller space, while working 
only with the pants graph.

The rank statement we obtain below is also obtained, for $S_{2,0}$, by
Brock-Masur \cite{BrockMasur:genus2}, as a consequence of an analysis of
the special properties of quasi-geodesics in the pants graph for
the genus 2 case. 

\bigskip

Recall that the \emph{Teichm\"{u}ller space} of a topological 
surface is the deformation 
space of finite area hyperbolic  structures which can be realized on 
that surface. Teichm\"{u}ller space has many natural metrics, here we 
consider the \emph{Weil-Petersson metric} which is a  K\"{a}hler 
metric with negative sectional curvature.

\begin{definition}{} The \emph{pants graph} of $S$ is a simplicial complex, 
    $\PP(S)$, with the following simplices:
    \begin{enumerate}
        \item {\bf Vertices}: one vertex for each pants decomposition 
	of $S$,  i.e., a top dimensional simplex in $\CC(S)$.
        \item {\bf Edges}: connect two pants decompositions 
	by an edge if they agree on all but one curve, and those 
	curves differ by an edge in the curve complex of the 
	complexity one subsurface (complementary to the rest of the 
	curves) in which they lie.
    \end{enumerate}    
\end{definition}

The following result of Brock \cite{Brock:wp} 
allows us to work with the pants graph in our study of Teichm\"{u}ller space.
\begin{theorem}{} $\PP(S)$ is quasi-isometric to the 
    Teichm\"{u}ller space of $S$ with the Weil-Petersson metric.
\end{theorem}

An important remark recorded in \cite{MasurMinsky:complex2} is that 
the pants graph is exactly what remains of the marking complex when 
annuli (and hence transverse curves) are ignored. Hence, one obtains 
the following version of Theorem~\ref{distance formula}:
\begin{theorem}{pants distance formula} 
    If $\mu,\nu\in \PP(S)$,  then there exists a constant $K(S)$, 
    depending only on the topological type of $S$, such that for each 
    $K>K(S)$ there exists $a\geq 1$ and $b\geq 0$ for which:
    $$
    d_{\PP(S)}(\mu,\nu) \approx_{a,b} 
    \sum_{\rm{non-annular \,}Y\subseteq S} 
    \Tsh K{d_{\CC(Y)}(\pi_{Y}(\mu),\pi_{Y}(\nu))}
    $$
\end{theorem}

\bigskip

We note that in \cite{Behrstock:asymptotic}, analogues of both 
Theorems~\ref{big implies small} and~\ref{F projection} are 
proven to hold for the pants graph 
of any surface of finite type. Further, by the above 
heuristic about ignoring annuli, one obtains product regions as 
produced for 
the mapping class group in Section~\ref{Qregions}. Again these
product regions  
are quasi-isometrically embedded 
with uniform constants; in the pants graph the identification is:
\begin{equation}\label{pants Q product}
Q_{\PP(S)}(\Delta) \homeo
\prod_{\rm{ non-annular \,}U\in \sigma(\Delta)} \PP(U).
\end{equation}
This identification leads to the main difference between the case of 
the pants graph and the mapping class group, namely, one obtains 
different counts of how many distinct factors occur on 
the right hand side of the above equation. In the mapping class group, 
this number is $3g+p-3$, whereas in the case of the pants graph, the 
count is easily verified to be $\lfloor\frac{3g+p-2}{2}\rfloor$.

As in the case of the mapping class group, one obtains:
\begin{lemma}{pants distance from Q}
If $\mu\in \PP(S)$ then 
$$
d(\mu,Q_{\PP(S)}(\Delta)) \approx 
\sum_{\substack{\tint{W}{\Delta}\\
W \rm{\, non-annular}
}} 
\Tsh K{d_{\CC(W)}(\mu,\Delta)}
$$
\end{lemma}

The remainder of the argument is completed as for the mapping class 
group, except for the count on the dimension of the separators. 
In the pants graph one obtains:
\begin{lemma}{pants separate points}
For any two points $x,y\in\AP$ there exists a closed set $L\subset
\AP$ which separates $x$ from $y$, and such that 
$\indhat(L) \le \lfloor\frac{3g+p-2}{2}\rfloor-1$. 
\end{lemma}

Thus, we have shown:

\state{Dimension theorem for Teichm\"{u}ller space.}
{Every locally-compact subset of an asymptotic cone of Teichm\"{u}ller 
space with the Weil-Petersson metric  
has topological dimension at most $\lfloor\frac{3g+p-2}{2}\rfloor$.}

The Rank Theorem for Teichmuller space now follows just as for the mapping
class group.

\bibliographystyle{ctmbib.bst}
\bibliography{/Users/jabehr/Documents/Math/TEX.macros.symbols/behrstock}

\end{document}